\documentclass[12pt]{amsart}
\usepackage{amsmath,amsthm,amssymb}
\usepackage[pdftex]{graphicx}
\usepackage{color}
\DeclareMathOperator\arctanh{arctanh}

\newtheorem{thm}{Theorem}[section]
\newtheorem{lem}[thm]{Lemma}
\newtheorem{prop}[thm]{Proposition}

\numberwithin{equation}{section}

\begin{document}

\begin{center}
 Morse index and symmetry-breaking bifurcation \\
 of positive solutions 
 to the one-dimensional Liouville type equation with a step function weight
\end{center}

\vspace{2ex}

\begin{center}
 Kanako Manabe \\[2ex]
 JG Corporation \\
 Konan 2--12--26, Minato-ku, Tokyo 108-0075, Japan \\[2ex]
 and \\[2ex]
 Satoshi Tanaka 
 \\[2ex]
 Mathematical Institute, Tohoku University \\ 
 Aoba 6--3, Aramaki, Aoba-ku, Sendai 980-8578, Japan \\
 satoshi.tanaka.d4@tohoku.ac.jp
\end{center}

\bigskip

{\bf Abstract.}
We consider the boundary value problem 
\begin{equation*}
 \left\{
  \begin{array}{l}
   u'' + \lambda h(x,\alpha) e^u = 0, \quad x \in (-1,1), \\[1ex]
   u(-1) = u(1) = 0,
  \end{array}
  \right.
\end{equation*}
where $\lambda>0$, $0<\alpha<1$, 
$h(x,\alpha)=0$ for $|x|<\alpha$, and 
$h(x,\alpha)=1$ for $\alpha \le |x| \le 1$.
We compute the Morse index of positive even solutions, and then
we prove the existence of an unbounded connected set of positive 
non-even solutions emanating from a symmetry-breaking bifurcation point.

\vspace{2ex}


\section{Introduction}

We consider the boundary value problem 
\begin{equation}
 \left\{
  \begin{array}{l}
   u'' +  \lambda h(x,\alpha) e^u = 0, \quad x \in (-1,1), \\[1ex]
    u(-1) = u(1) = 0,
    \label{P}
  \end{array}
  \right.
\end{equation}
where $\lambda>0$, $0<\alpha<1$, $h(x,\alpha)=0$ for $|x|<\alpha$, and 
$h(x,\alpha)=1$ for $\alpha \le |x| \le 1$.
It is said to be a {\it solution} of \eqref{P} if
\begin{equation*}
 u \in C^2([-1,1]\setminus\{-\alpha,\alpha\}) \cap C^1[-1,1], \quad
 u(-1)=u(1)=0,
\end{equation*}
and it satisfies 
\begin{equation*}
 u'' +  \lambda h(x,\alpha) e^u = 0, 
 \quad x \in (-1,1)\setminus\{-\alpha,\alpha\}.
\end{equation*}

Joseph and Lundgren \cite{JL} elucidated the complete structure of 
positive radial solutions to the Liouville equation in the unit ball
\begin{equation}
 \left\{
  \begin{array}{cl}
   \Delta u + \lambda e^u = 0 & \mbox{in}\ \Omega, \\[1ex]
    u = 0 & \mbox{on}\ \partial \Omega,
  \end{array}
  \right.
 \label{Liouville}
\end{equation}
where $\Omega=B:=\{ x \in \mathbb{R}^N : |x|<1 \}$ and $N\ge 1$.
By a celebrated theorem by Gidas, Ni and Nirenberg \cite{GNN}, every positive 
solution of \eqref{Liouville} is radially symmetric when $\Omega=B$.
However, when $\Omega$ is an annulus 
$A:=\{ x \in \mathbb{R}^N : a < |x| < b \}$, $0<a<b$, 
the Gidas-Ni-Nirenberg theorem can not be applied, and 
problem \eqref{Liouville} may have non-radial solutions.
In particular, Lin \cite{Lin} found infinitely many symmetry-breaking 
bifurcation points on the branch of positive radial solutions of 
\eqref{Liouville} with $\Omega=A$ and $N=2$.
We also can not apply the Gidas-Ni-Nirenberg theorem to 
the following problem
\begin{equation}
 \left\{
  \begin{array}{cl}
   \Delta u + \lambda |x|^l e^u = 0 & \mbox{in}\ B, \\[1ex]
    u = 0 & \mbox{on}\ \partial B,
  \end{array}
  \right.
 \label{LT} 
\end{equation}
where $\lambda>0$, $l>0$ and $N \ge 1$.
Indeed, Miyamoto \cite{Miya} presented the complete result of structure to 
\eqref{LT} when $N=2$, and as a sequel to it, we conclude that 
the number of symmetry-breaking bifurcation points goes to infinity as 
$l\to\infty$.
For the case $N=1$, that is, for the problem
\begin{equation}
 \left\{
  \begin{array}{l}
   u'' + \lambda |x|^l e^u = 0, \quad x \in (-1,1), \\[1ex]
    u(-1) = u(1) = 0,
  \end{array}
  \right.
   \label{N=1}
\end{equation}
at least one symmetry-breaking bifurcation point is found in \cite{T2017}. 
However, there is an open part in the result of \cite[Conjecture 1.1]{T2017}.
Specifically, it is conjectured that such a symmetry-breaking bifurcation 
point is unique, but this remains unknown.
Therefore, instead of \eqref{N=1}, we consider problem \eqref{P} in this paper.
Problem \eqref{N=1} can be considered as an approximation of problem \eqref{P}
in the following sense.
For every $q\ge 1$, we have
\begin{align*}
 \int_{-1}^1 ||x|^l-h(x,e^{1/l})|^q dx 
  = \frac{2e^{q+\frac{1}{l}}}{lq+1} 
  +2\int_{e^{1/l}}^1 (1-x^l)^q dx \to 0
\end{align*}
as $l \to \infty$. 
Since the differential equation in \eqref{P} is piece-wise autonomous,
we expect to obtain detailed results.
In this paper, we will demonstrate the existence and the uniqueness of 
symmetry-breaking bifurcation point of \eqref{P}.

There has been much current interest in studying
boundary value problem with a step function weight.
Kajikiya \cite{Kaj2021, Kaj2022JMAA, Kaj2022JMSJ, Kaj2023} has thoroughly 
studied the Moore-Nehari equation
\begin{equation*}
 u'' + h(x,\alpha) |u|^{p-1}u = 0, \ \ x \in (-1,1); \quad
 u(-1) = u(1) = 0,
\end{equation*}
where $h$ is the same function as in \eqref{P}, $p>0$ and $p\ne1$.
See also \cite{KST}.
Cubillos, L\'{o}pez-G\'{o}mez and Tellini \cite{CLGT2022, CLGT2023, CLGT2024}
showed that the solution structure of the problem
\begin{equation*}
 -u'' = \lambda u - a(x) u^3, \quad x\in(0,1); \quad u(0)=u(1)=0
\end{equation*}
is very interesting, where $\lambda\in\mathbb{R}$ is a parameter and $a$ is a 
piece-wise continuous function.
See also, L\'{o}pez-G\'{o}mez, Mu\~{n}oz-Hern\'{a}ndez and Zanolin \cite{LGMHZ},
 L\'{o}pez-G\'{o}mez and Rabinowitz \cite{LGR2017, LGR2020}, and Kan \cite{Kan}.

Let $w$ be the unique solution of the initial value problem
\begin{equation*}
 \left\{
  \begin{array}{l}
   w'' + e^w = 0, \quad x \in \mathbb{R}, \\[1ex]
   w(0)=w'(0)=0,
  \end{array}
 \right.
\end{equation*}
that is,
\begin{equation}
 w(x)= \log \left( 1 - \tanh^2 \frac{x}{\sqrt{2}} \right).
 \label{w}
\end{equation}
Then $w$ is an even function and satisfies
\begin{equation*}
 w(x)<0, \ w'(x)<0, \ w''(x)<0 \quad \mbox{for}\ x>0 
\end{equation*}
and $\lim_{x\to\infty} w(x)=-\infty$.
Let $\eta$ be the inverse function of $\beta=-w(x)$ on $[0,\infty)$.
Namely,
\begin{equation*}
 \eta(\beta)= \sqrt{2} \arctanh \sqrt{1-e^{-\beta}}.
\end{equation*}
Then $\eta \in C^2(0,\infty)$, $\eta(\beta)>0$,
$\eta'(\beta)>0$ for $\beta>0$, $\eta(0)=0$, and
$\lim_{\beta\to\infty} \eta(\beta)=\infty$.
We set
\begin{equation}
 \Lambda(\beta)=(1-\alpha)^{-2}e^{-\beta}[\eta(\beta)]^2
  \label{Kormanlam}
\end{equation}
and
\begin{equation}
 U(x;\beta)=
 \left\{
  \begin{array}{ll}
   \beta & |x|<\alpha, \\[1ex]
   w\left(\eta(\beta)\dfrac{|x|-\alpha}{1-\alpha}\right)+\beta, 
   & \alpha \le |x| \le 1.
  \end{array}
 \right. 
 \label{Kormansol}
\end{equation}
Since the function $f(x):=x\arctanh x$ is strictly increasing in $x\in[0,1)$ 
and satisfies $f(0)=0$ and $\lim_{x\to1} f(x)=\infty$, the equation
\begin{equation}
 2\eta'(\beta)-\eta(\beta) = 0
  \label{2eta'=eta}
\end{equation}
that is,
\begin{equation}
 \sqrt{1-e^{-\beta}} \arctanh \sqrt{1-e^{-\beta}} = 1
  \label{betaast}
\end{equation}
has the unique positive root.
Let $\beta_1>0$ be the unique positive root of \eqref{betaast}.
Then $\beta_1=1.18...$ and 
$\Lambda(\beta)$ satisfies $\Lambda'(\beta)>0$ for $0<\beta<\beta_1$, 
$\Lambda'(\beta_1)=0$, $\Lambda'(\beta)<0$ for $\beta>\beta_1$, and
\begin{equation}
    \lim_{\beta\to0^+} \Lambda(\beta) 
  = \lim_{\beta\to\infty} \Lambda(\beta) = 0.
   \label{limlambda}
\end{equation}
It is easy to prove the following result.

\begin{lem}\label{Uisevensol}
 For each $\beta>0$, $U(x;\beta)$ is a positive even solution of
 \eqref{P} with $\lambda=\Lambda(\beta)$ and $\|U\|_\infty=\beta$.
\end{lem}

Here and hereafter we use the notation: 
$||u||_\infty=\sup_{x \in [-1,1]}|u(x)|$.

Moreover, we have the following result, which will be shown in Section 4.

\begin{prop}\label{evensolutions}
 For each $\beta>0$, there exists the unique $(\lambda,u)$ such that 
 $\lambda>0$ and $u$ is a positive even solution of \eqref{P} with 
 $||u||_\infty=\beta$, that is, 
 $(\lambda,u)=(\Lambda(\beta),U(\,\cdot\,;\beta))$.
 Moreover, the following \textup{(i)--(iii)} hold\textup{:}
 \begin{enumerate}
  \item if $0<\lambda<\lambda_1$, then \eqref{P} has exactly two 
	positive even solutions\textup{;}
  \item if $\lambda=\lambda_1$, then \eqref{P} has the unique 
	positive even solution\textup{;}
  \item if $\lambda>\lambda_1$, then \eqref{P} has no positive 
	 even solution,	
 \end{enumerate}
 where $\lambda_1=\Lambda(\beta_1)$ and $\beta_1$ is the unique positive root 
 of \eqref{betaast}.
\end{prop}

Let $m(\beta)$ be the Morse index of $U(x;\beta)$, that is, 
the number of negative eigenvalues $\mu$ to
\begin{equation}
 \left\{
 \begin{array}{l}
  \phi'' + \Lambda(\beta) h(x,\alpha) e^{U(x;\beta)} \phi + \mu \phi = 0,
   \quad x \in (-1,1), \\[1ex]
  \phi(-1) = \phi(1) = 0.
 \end{array}
	\right.
  \label{Morse}
\end{equation}
A solution $U(x;\beta)$ is said to be degenerate if
$\mu=0$ is an eigenvalue of \eqref{Morse}.
Otherwise, it is said to be nondegenerate.

We denote by $\mu_k(\beta)$ the $k$-th eigenvalue of \eqref{Morse}.
We recall that 
\[
 \mu_1(\beta) < \mu_2(\beta) < \cdots < \mu_k(\beta) <
 \mu_{k+1}(\beta) < \cdots, \quad
 \lim_{k\to\infty} \mu_k(\beta) = \infty,
\]
problem \eqref{Morse} has 
no other eigenvalues, an eigenfunction $\phi_k$ corresponding to
$\mu_k(\beta)$ is unique up to a constant, and $\phi_k$ has exactly
$k-1$ zeros in $(-1,1)$.
We find that $\mu_k \in C(0,\infty)$. (See, for example, \cite{Kos}.)

The following theorem is the main result of this paper.

\begin{thm}\label{main}
 Let $(\Lambda(\beta),U(x;\beta))$ be as in 
 \eqref{Kormanlam}--\eqref{Kormansol} and 
 let $\beta_1$ be the unique positive root of \eqref{betaast}.
 Then there exists a constant $\beta_2$ such that 
 $\beta_2>\beta_1$ and the following \textup{(i)--(v)} hold\textup{:}
 \begin{enumerate}
  \item if $0<\beta<\beta_1$, then $m(\beta)=0$ and $U(x;\beta)$ is
	nondegenerate\textup{;}
  \item if $\beta=\beta_1$, then $m(\beta)=0$ and $U(x;\beta)$ is
	degenerate\textup{;}
  \item if $\beta_1<\beta<\beta_2$, then $m(\beta)=1$ and $U(x;\beta)$
	is nondegenerate\textup{;}
  \item if $\beta=\beta_2$, then $m(\beta)=1$ and $U(x;\beta)$ is
	degenerate\textup{;}
  \item if $\beta>\beta_2$, then $m(\beta)=2$ and $U(x;\beta)$ is
	nondegenerate.
 \end{enumerate}
 Moreover, there exists an unbounded connected set 
 $\mathcal{C}\subset (0,\infty)\times C^1_0[-1,1]$ which emanates from
 $(\Lambda(\beta_2),U(\,\cdot\,;\beta_2))$.
 For each point $(\lambda,u) \in \mathcal{C}$,
 $u$ is a positive non-even solution of \eqref{P} and satisfies
 \begin{equation*}
  \Lambda(||u||_\infty) \le \lambda < 4 \Lambda(||u||_\infty),
 \end{equation*}
 and there exist a sequence $\{(\lambda_n,u_n)\}\subset\mathcal{C}$ such that
 $\lambda_n\to0$ and $||u_n||_\infty\to\infty$.
\end{thm}

\begin{center}
 \includegraphics[width=5.5cm]{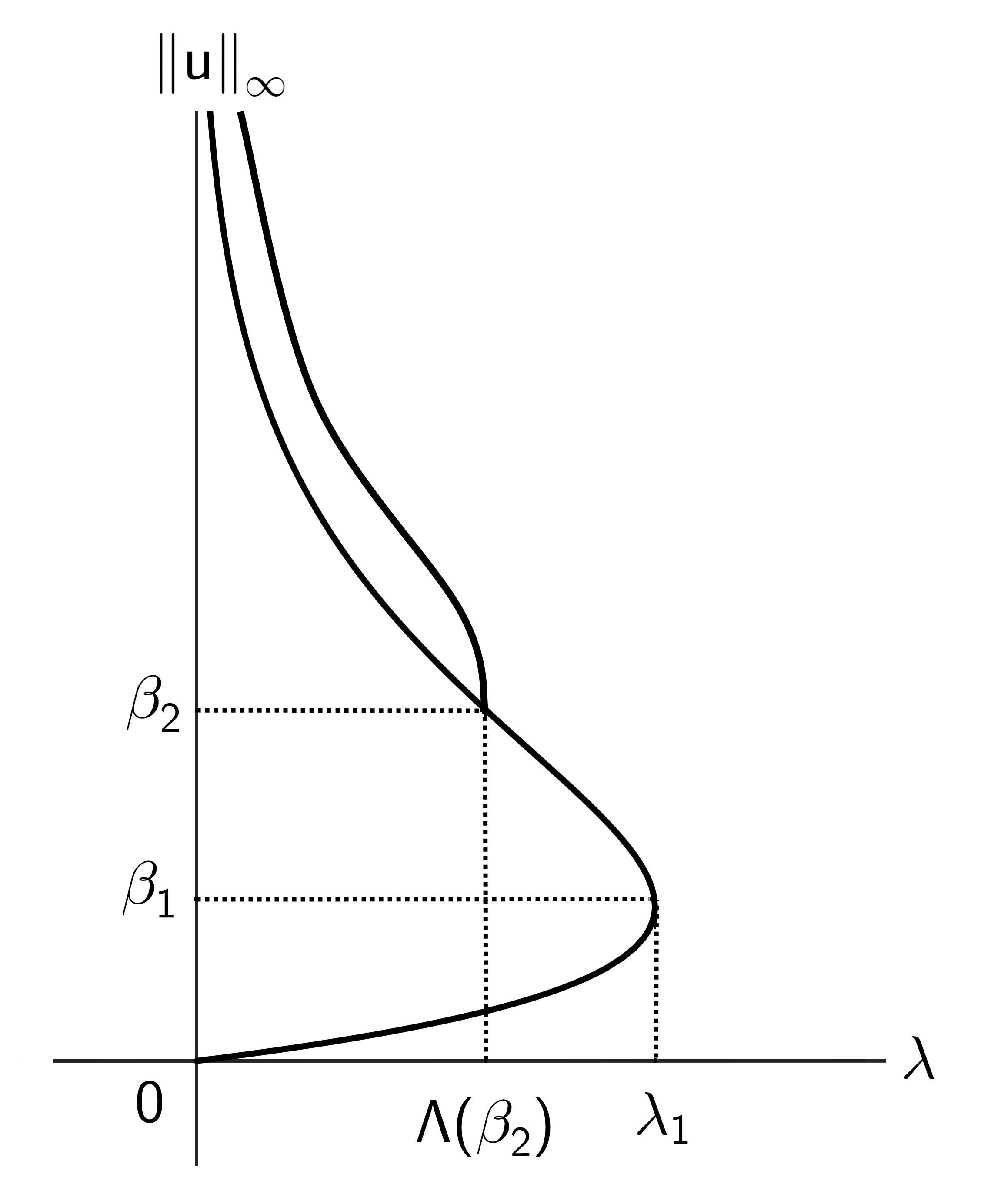} \\
 Illustration of Theorem \ref{main}
\end{center}

\medskip

In Section 2, we prove Proposition \ref{evensolutions}.
We study the eigenvalues $\mu_1(\beta)$, $\mu_3(\beta)$ and $\mu_2(\beta)$ 
in Sections 3 and 4, respectively.
In Section 5, we give a proof of Theorem \ref{main}.

\section{Proof of Proposition \ref{evensolutions}}

In this section, we provide the proof of Proposition \ref{evensolutions}.

\begin{proof}[Proof of Proposition \ref{evensolutions}]
 Let $\beta>0$ be fixed.
 Let $u$ be a positive even solution of \eqref{P} with $||u||_\infty=\beta$. 
 Then $u$ and $\lambda$ must satisfy
 \begin{gather*}
  u'' + \lambda e^u= 0, \ x\in(-1,-\alpha); \quad 
  u(-\alpha)=\beta, \ u'(-\alpha)=0; \\
  u(x)=\beta, \ x \in (-\alpha,\alpha) \\
  u'' + \lambda e^u= 0, \ x\in(\alpha,1); 
  \quad u(\alpha)=\beta, \ u'(\alpha)=0; \\
  u(-1)=u(1)=0,
 \end{gather*}
 that is,
 \begin{equation*}
  u(x) =
   \left\{
  \begin{array}{ll}
   \beta & |x|<\alpha, \\[1ex]
   w\left(\sqrt{\lambda e^\beta}(|x|-\alpha)\right) + \beta
   & \alpha \le |x| \le 1.
  \end{array}
 \right. 
 \end{equation*}
 where $w$ is the function defined by \eqref{w}.
 Since $u(-1)=u(1)=0$, we find that $\lambda$ satisfies
 \begin{equation*}
  -w\left(\sqrt{\lambda e^\beta}(1-\alpha)\right) = \beta.
 \end{equation*}
 Recalling that $\eta$ is the inverse function of $\beta=-w(x)$ on $[0,\infty)$,
 we have
 \begin{equation*}
  \sqrt{\lambda e^\beta}(1-\alpha) = \eta(\beta),
 \end{equation*}
 which means that $\lambda$ must be $\Lambda(\beta)$.
 Consequently, there exists the unique $(\lambda,u)$ such that 
 $\lambda>0$ and $u$ is a positive even solution of \eqref{P} with 
 $||u||_\infty=\beta$, that is, 
 $(\lambda,u)=(\Lambda(\beta),U(\,\cdot\,;\beta))$.

 We recall that $\Lambda(\beta)$ satisfies $\Lambda'(\beta)>0$ for 
 $0<\beta<\beta_1$, $\Lambda'(\beta_1)=0$, $\Lambda'(\beta)<0$ for 
 $\beta>\beta_1$, and \eqref{limlambda} holds.
 If $0<\lambda<\lambda_1=\Lambda(\beta_1)$, then there exist 
 exactly two positive numbers $\beta_*$ and $\beta^*$ such that 
 $\Lambda(\beta_*)=\Lambda(\beta^*)=\lambda$, and hence
 $U(x;\beta_*)$ and $U(x;\beta^*)$ are positive even solutions of 
 \eqref{P} and there is no other positive even 
 solution.
 If $\lambda=\lambda_1$, then 
 $U(x;\beta_1)$ is the unique positive even solution of \eqref{P}.
 If $\lambda>\lambda_1$, then there is no $\beta>0$ such that 
 $\Lambda(\beta)=\lambda$, which means that \eqref{P} has no 
 positive even solution.
\end{proof}

\section{The first and third eigenvalues} 

In this section, we study the eigenvalues $\mu_1(\beta)$ and $\mu_3(\beta)$
of the linearized problem \eqref{Morse}.
We define 
\begin{equation*}
 \psi(x;\beta) := (x-\alpha) U'(x;\beta) + 2.
\end{equation*}
We easily check that the next lemma follows from Lemma \ref{Uisevensol}.

\begin{lem}\label{psi1}
 The function $\psi(x;\beta)$ is a solution of 
  \begin{equation}
   \psi'' + \Lambda(\beta) h(x,\alpha) e^{U(x;\beta)} \psi = 0.
    \label{psi}
  \end{equation}
\end{lem}

\begin{lem}\label{psi2} 
 Let $\beta_1>0$ be the unique positive root of \eqref{betaast}.
 Then $\psi(x;\beta)$ satisfies the following \textup{(i)--(iii):}
 \begin{enumerate}
  \item if $0<\beta<\beta_1$, then $\psi(x;\beta)>0$ for
	$x \in [-1,1]$\textup{;} 
  \item $\psi(x;\beta_1)>0$ for $x \in (-1,1)$ and
	$\psi(-1;\beta_1)=\psi(1;\beta_1)=0$\textup{;} 
  \item if $\beta>\beta_1$, then $\psi(x;\beta)$ has exactly two
	zeros in $(-1,1)$, \linebreak
	$\psi(-1;\beta)<0$ and $\psi(1;\beta)<0$. 
 \end{enumerate} 
\end{lem}

\begin{proof}
 Since
 \begin{equation*}
  \psi(x;\beta) = \eta(\beta)\dfrac{x-\alpha}{1-\alpha} 
     w'\left( \eta(\beta)\dfrac{x-\alpha}{1-\alpha} \right)+2, 
     \quad x \in [\alpha,1],
 \end{equation*}
 and 
 \begin{equation}
  (xw'(x))' = w'(x) + x w''(x) = w'(x) - xe^{w(x)} < 0,
  \quad x > 0,
   \label{(xw')'<0} 
 \end{equation}
 we conclude that $\psi(x;\beta)$ is strictly decreasing in $x \in [\alpha,1]$.
 We note that $\psi(x;\beta)=2$ for $x \in (-\alpha,\alpha)$ and $\psi(x;\beta)$ 
 is an even function, which implies
 \begin{equation*}
  \min_{x \in [-1,1]} \psi(x;\beta) = \psi(1;\beta)
  = \eta(\beta) w'(\eta(\beta))+2.
 \end{equation*}
 Recalling $\eta(\beta)$ is the inverse function of $\beta=-w(x)$, 
 we have
 \begin{equation*}
  \eta(\beta) w'(\eta(\beta))+2
  = -\frac{\eta(\beta)}{\eta'(\beta)} + 2 
 \end{equation*}
 Since $\beta_1$ is a solution of \eqref{2eta'=eta}, we get 
 \begin{equation*}
  \eta(\beta_1) w'(\eta(\beta_1))+2 = 0.
 \end{equation*}
 From \eqref{(xw')'<0}, it follows that $\eta(\beta) w'(\eta(\beta))+2$ 
 is strictly decreasing in $\beta>0$.
 Therefore, Lemma \ref{psi2} holds.
\end{proof}

We will employ the following version of the Sturm comparison theorem.
See for example Kajikiya \cite[Lemma 2.1]{Kaj2021}.

\begin{lem}\label{Sturm}
 Suppose that $h$, $H \in L^1 (a,b)$ and $h(x) \le H(x)$ a.e. in $(a,b)$. 
 Let $u$, $v \in W^{2,1}(a,b)$ be non-trivial solutions of 
 \begin{equation*}
  u'' + h(x) u =0
 \end{equation*}
 and
 \begin{equation*}
  v'' + H(x) v =0
 \end{equation*}
 in $(a,b)$, respectively.
 Suppose that $u(a)=u(b)=0$ and $u(x)\ne0$ in $(a,b)$.
 Then one of the following alternatives holds\textup{:}
 \begin{enumerate}
  \item $v$ has at least one zero in $(a,b)$\textup{;}
  \item $v$ is a constant multiple of $u$.
 \end{enumerate}
 The second alternative implies that $h(x)=H(x)$ a.e. in $(a,b)$. 
 Therefore, if the set of points $x$ satisfying $h(x)<H(x)$ has 
 a positive Lebesgue measure, then the assertion \textup{(i)} only holds.
\end{lem}

\begin{lem}\label{mu1}
 Let $\beta_1>0$ be the unique positive root of \eqref{betaast}.
 Then the first eigenvalue $\mu_1(\beta)$ of \eqref{Morse} satisfies the
 following \textup{(i)--(iii):}
 \begin{enumerate}
  \item $\mu_1(\beta)>0$ for $0<\beta<\beta_1$\textup{;} 
  \item $\mu_1(\beta_1)=0$\textup{;} 
  \item $\mu_1(\beta)<0$ for $\beta>\beta_1$\textup{.} 
 \end{enumerate} 
\end{lem}

\begin{proof}
 Let $\phi_1$ be an eigenfunction corresponding to $\mu_1(\beta)$.
 Then $\phi_1(x) \ne 0$ on $(-1,1)$ and $\phi_1(-1)=\phi_1(1)=0$.

 (i)
 Assume to the contrary that $\mu_1(\beta_0)\le 0$ for some 
 $\beta_0 \in (0,\beta_1)$.
 From Lemma \ref{Sturm}, it follows that every solution of \eqref{psi} with
 $\beta=\beta_0$ has at least one zero in $[-1,1]$,
 which contradicts (i) of Lemma \ref{psi2}.
 Therefore, $\mu_1(\beta)>0$ for $0<\beta<\beta_1$.

 (ii) 
 By (ii) of Lemma \ref{psi2}, we conclude that $\psi(x;\beta_1)$
 is an eigenfunction corresponding to $\mu_1(\beta_1)$, which means
 $\mu_1(\beta_1)=0$.

 (iii)
 Suppose that $\mu_1(\beta_3) \ge 0$ for some $\beta_3>\beta_1$.
 By (iii) of Lemma \ref{psi2}, there exists $x_0 \in (0,1)$ such that
 $\psi(-x_0;\beta_3)=\psi(x_0;\beta_3)=0$.
 Lemma \ref{Sturm} implies that every solution of
 \begin{equation*}
  \phi''
   + [\Lambda(\beta_3) h(x,\alpha) e^{U(x;\beta_3)} + \mu_1(\beta_3)] \phi = 0
 \end{equation*}
 has at least one zero in $[-x_0,x_0]$.
 This contradicts the fact that $\phi_1(x) \ne 0$ on $(-1,1)$. 
 Hence, $\mu_1(\beta)<0$ for $\beta>\beta_1$.
\end{proof}

\begin{lem}\label{mu3}
 The third eigenvalue $\mu_3(\beta)$ of \eqref{Morse} is positive for
 $\beta>0$.
\end{lem}

\begin{proof}
 Suppose that $\mu_3(\beta_0) \le 0$ for some $\beta_0>0$.
 Let $\phi_3$ be an eigenfunction of \eqref{Morse} corresponding to
 $\mu_3(\beta_0)$.
 Then $\phi_3(-1)=\phi_3(1)=0$ and $\phi_3$ has exactly two zeros in
 $(-1,1)$.
 Lemma \ref{Sturm} implies that every solution of \eqref{psi} with $\beta=\beta_0$
 has at least three zeros in $[-1,1]$.
 On the other hand, from Lemmas \ref{psi1} and \ref{psi2}, it follows that
 $\psi(x;\beta_0)$ is a solution of \eqref{psi} and 
 has at most two zeros in $[-1,1]$.
 This is a contradiction.
 Consequently, $\mu_3(\beta)>0$ for $\beta>0$.
\end{proof}

\section{The second eigenvalue}

For each $\beta>0$, we define the function $\varphi(x;\beta)$ by
\begin{equation*}
 \varphi(x;\beta) = \left\{
  \begin{array}{ll}
   (c_1(\beta)-\alpha x)U'(x;\beta) - 2\alpha, & -1 \le x \le -\alpha, \\[1ex]
   2x, & |x|<\alpha, \\[1ex]
   (c_1(\beta)+\alpha x)U'(x;\beta) + 2\alpha, & \alpha \le x \le 1,
  \end{array}
 \right.
\end{equation*}
where
\begin{equation*}
 c_1(\beta) := -\frac{2}{e^\beta\Lambda(\beta)}-\alpha^2 
      = -\frac{2(1-\alpha)^2}{[\eta(\beta)]^2}-\alpha^2. 
\end{equation*}
Then $\varphi(\,\cdot\,;\beta)\in C^1[-1,1]\cap C^2([-1,1]\setminus\{-\alpha,\alpha\})$ 
and is a solution of \eqref{psi} and an odd function for each fixed $\beta>0$.

\begin{lem}\label{phi1}
 Let $\beta>0$.
 The function $\varphi(x;\beta)$ has at most one zero in $(0,1]$.
\end{lem}

\begin{proof}
 Assume that $\varphi(x;\beta)$ has two zeros in $(0,1]$.
 We note that $x=0$ is also a zero of $\varphi(x;\beta)$.
 Lemmas \ref{psi1} and \ref{Sturm} imply that $\psi$ has two zeros in $(0,1]$.
 Since $\psi$ is an even function, we see that $\psi$ has at least four 
 zeros in $[-1,1]$.
 This contradicts Lemma \ref{psi2}.
\end{proof}

\begin{lem}\label{phi2}
 There exists $\beta_2>0$ such that 
 $\varphi(1;\beta)>0$ for $0<\beta<\beta_2$,
 $\varphi(1;\beta_2)=0$ and
 $\varphi(1;\beta)<0$ for $\beta>\beta_2$.
\end{lem}

\begin{proof}
 We observe that
 \begin{align*}
  \varphi(1;\beta) & = (c_1(\beta)+\alpha) U'(1;\beta) + 2\alpha \\
  & = \left( -\frac{2(1-\alpha)^2}{[\eta(\beta)]^2}-\alpha^2 + \alpha \right)
      w'(\eta(\beta))\frac{\eta(\beta)}{1-\alpha} + 2\alpha \\
  & = -\left( -\frac{2(1-\alpha)^2}{[\eta(\beta)]^2}-\alpha^2 + \alpha \right)
      \frac{\sqrt{2}\eta(\beta)}{1-\alpha} 
           \tanh \frac{\eta(\beta)}{\sqrt{2}} + 2\alpha \\
  & = g\left( \frac{\eta(\beta)}{\sqrt{2}} \right),
 \end{align*}
 where 
 \begin{equation*}
  g(x):= 2 ( 1-\alpha-\alpha x^2 ) \frac{1}{x} \tanh x + 2\alpha.
 \end{equation*}
 We observe that $\lim_{x\to0^+} g(x)=2$ and $\lim_{x\to\infty} g(x)=-\infty$.
 Moreover, $g'(x)<0$ for $x>0$.
 Indeed, 
 \begin{equation*}
  \frac{x^2}{2} g'(x)
   = -(1-\alpha+\alpha x^2) \tanh x + (1-\alpha-\alpha x^2) x(1-\tanh^2 x)
 \end{equation*}
 and using the inequality 
 \begin{equation*}
  \tanh x > x(1-\tanh^2 x), \quad x>0,   
 \end{equation*}
 we have
 \begin{align*}
  \frac{x^2}{2} g'(x)
  & < -(1-\alpha+\alpha x^2) x(1-\tanh^2 x) 
      +(1-\alpha-\alpha x^2) x(1-\tanh^2 x) \\
  & = -2\alpha x^3 (1-\tanh^2 x), \quad x>0.
 \end{align*}
 Recalling $\eta(\beta)>0$, $\eta'(\beta)>0$ for $\beta>0$, $\eta(0)=0$, and
 $\lim_{\beta\to\infty} \eta(\beta)=\infty$, we conclude that
 $\lim_{\beta\to0}\varphi(1;\beta)=2$,
 $\lim_{\beta\to\infty}\varphi(1;\beta)=-\infty$ and 
 and hence there exists $\beta_2>0$ such that 
 $\varphi(1;\beta)>0$ for $0<\beta<\beta_2$,
 $\varphi(1;\beta_1)=0$ and
 $\varphi(1;\beta)<0$ for $\beta>\beta_2$.
\end{proof}

Since $\varphi(x;\beta)=2x>0$ for $0<x\le\alpha$, combining Lemmas \ref{phi1}
and \ref{phi2}, we obtain the following lemma.

\begin{lem}\label{phi3}
 Let $\beta_2$ be as in Lemma \ref{phi2}.
 Then $\varphi(x;\beta)$ satisfies the following \textup{(i)--(iii):}
 \begin{enumerate}
  \item if $0<\beta<\beta_2$, then $\varphi(x;\beta)>0$ for
	$x \in (0,1]$\textup{;} 
  \item $\varphi(x;\beta_2)>0$ for $x \in (0,1)$ and
	$\varphi(1;\beta_2)=0$\textup{;} 
  \item if $\beta>\beta_2$, then $\varphi(x;\beta)$ has the unique zero in 
	$(0,1)$ and $\varphi(1;\beta)<0$. 
 \end{enumerate} 
\end{lem}

\begin{lem}\label{mu2}
 Let $\beta_1>0$ be the unique positive root of \eqref{betaast} and
 let $\beta_2$ be as in Lemma \ref{phi2}.
 Then the second eigenvalue $\mu_2(\beta)$ of \eqref{Morse} satisfies the
 following \textup{(i)--(iii):}
 \begin{enumerate}
  \item $\mu_2(\beta)>0$ for $0<\beta<\beta_2$\textup{;} 
  \item $\mu_2(\beta_2)=0$\textup{;} 
  \item $\mu_2(\beta)<0$ for $\beta>\beta_2$\textup{.} 
 \end{enumerate}
 Moreover, $0<\beta_1<\beta_2$.
\end{lem}

\begin{proof}
 The proof is similar to the proof of Lemma \ref{mu1}.
 
 Let $\phi_2$ be an eigenfunction corresponding to $\mu_2(\beta)$.
 Then $\phi_2(x)$ has the unique zero in $(-1,1)$ and $\phi_1(-1)=\phi_1(1)=0$.
 We recall that $\varphi(x;\beta)$ is an odd function. 

 (i)
 Suppose that $\mu_2(\beta_0)\le 0$ for some $\beta_0 \in (0,\beta_2)$.
 From Lemma \ref{Sturm}, it follows that every solution of \eqref{psi} with
 $\beta=\beta_0$ has at least two zeros in $[-1,1]$.
 This contradicts (i) of Lemma \ref{phi3}.
 Therefore, $\mu_2(\beta)>0$ for $0<\beta<\beta_2$.

 (ii) 
 From (ii) of Lemma \ref{phi3}, it follows that $\varphi(x;\beta_2)$
 is an eigenfunction corresponding to $\mu_2(\beta_2)$ and hence 
 $\mu_2(\beta_2)=0$.

 (iii)
 Assume that $\mu_2(\beta_3) \ge 0$ for some $\beta_3>\beta_2$.
 By (iii) of Lemma \ref{phi3}, there exists $x_0 \in (0,1)$ such that
 $\varphi(-x_0;\beta_3)=\varphi(x_0;\beta_3)=0$.
 We note that $x=0$ is also a zero of $\varphi(x;\beta_3)=0$.
 By Lemma \ref{Sturm}, we see that every solution of
 \begin{equation*}
  \phi''
   + [\Lambda(\beta_3) h(x,\alpha) e^{U(x;\beta_3)} + \mu_2(\beta_3)] \phi = 0  
 \end{equation*}
 has at least two zeros in $[-x_0,x_0]$.
 This contradicts the fact that $\phi_2(x)$ has the unique zero in $(-1,1)$.
 We conclude that $\mu_2(\beta)<0$ for $\beta>\beta_2$.

 Since $\mu_1(\beta)<\mu_2(\beta)$ for $\beta>0$, we have
 $0=\mu_1(\beta_1)<\mu_2(\beta_1)$.
 Then $\beta_1$ must satisfy $0<\beta_1<\beta_2$.
\end{proof}

\section{Proof of the main result}

In this section, we prove Theorem \ref{main}.
From Lemmas \ref{mu1}, \ref{mu3} and \ref{mu2}, (i)--(v) of Theorem \ref{main} 
follows immediately.
To prove the remaining parts of Theorem \ref{main}, 
we need the following lemmas

\begin{lem}\label{lemlambda<}
 Let $u$ be a positive solution of \eqref{P}.
 Then 
 \begin{equation*}
  \Lambda(||u||_\infty) \le \lambda < 4 \Lambda(||u||_\infty),
 \end{equation*}
 where $\Lambda(\beta)$ is the function defined by \eqref{Kormanlam}.
\end{lem}

\begin{proof}
 Since $u''(x)=0$ for $-\alpha<x<\alpha$, there exists 
 $m \in (-1,\alpha]\cup[\alpha,1)$ such that $u'(m)=0$ and 
 $u(m)=||u||_\infty$.
 Without loss of generality, we suppose that $m\in[\alpha,1)$,
 since $u(-x)$ is also a positive solution of \eqref{P}.
 It follows that 
 \begin{equation*}
  v(x):= w\Bigl(\sqrt{\lambda e^{||u||_\infty}}(x-m)\Bigr)+||u||_\infty
 \end{equation*}
 is a solution of
 \begin{equation*}
  v''+\lambda e^v =0, \quad \alpha<x\le 1, \quad
  v(m)=||u||_\infty, \ v'(m)=0.
 \end{equation*}
 By the uniqueness of initial value problems, we find that
 $u(x)\equiv v(x)$ for $x \in [\alpha,1]$.
 In particular, $v(1)=u(1)=0$.
 Therefore,
 \begin{equation*}
  -w\Bigl(\sqrt{\lambda e^{||u||_\infty}}(1-m)\Bigr)=||u||_\infty,
 \end{equation*}
 which is equivalent to
 \begin{equation*}
  \sqrt{\lambda e^{||u||_\infty}}(1-m) = \eta(||u||_\infty), 
 \end{equation*}
 that is,
 \begin{equation}
  \lambda = (1-m)^{-2} e^{-||u||_\infty} [\eta(||u||_\infty)]^2
   = \left(\frac{1-\alpha}{1-m}\right)^2 \Lambda(||u||_\infty).
   \label{lambda=}
 \end{equation}
 Since $m \ge \alpha$, we have $\lambda\ge\Lambda(||u||_\infty)$.
 Recalling $w$ is an even function, we see that $v(x)=v(2m-x)$, which implies 
 $v(2m-1)=v(1)=0$.
 By $m<1$, we have $2m-1<1$.
 Since $v(x)\equiv u(x)>0$ for $x\in[\alpha,1)$, we conclude that $2m-1<\alpha$.
 By \eqref{lambda=}, we get $\lambda < 4 \Lambda(||u||_\infty)$.
\end{proof}

\begin{lem}\label{||u||<1}
 Let $\lambda>0$.
 Then problem \eqref{P} with $||u||_\infty \le 1$ has at most one
 positive solution.
\end{lem}

\begin{proof}
 Assume that problem \eqref{P} has two distinct positive solutions 
 $u$ and $v$ for which $||u||_\infty \le 1$ and $||v||_\infty \le 1$.
 Then there exist $x_1$, $x_2 \in [-1,1]$ such that $u(x_1)=v(x_1)$,
 $u(x_2)=v(x_2)$ and $u(x)\ne v(x)$ for $x_1<x<x_2$.
 We suppose that $u(x)>v(x)$ for $x_1<x<x_2$.
 Then $u'(x_1)>v'(x_1)$ and $u'(x_2)<v'(x_2)$.
 We find that $(x_1,x_2)\not\subset(-\alpha,\alpha)$.
 Indeed, if $(x_1,x_2)\subset(-\alpha,\alpha)$, then $u''(x)=v''(x)=0$
 for $x_1<x<x_2$, and hence $u(x)=u'(x_1)(x-x_1)+u(x_1)$ and
 $v(x)=v'(x_1)(x-x_1)+v(x_1)$.
 Therefore,
 \begin{equation*}
  u(x_2) = u'(x_1)(x_2-x_1)+u(x_1) > v'(x_1)(x_2-x_1)+v(x_1) = v(x_2),
 \end{equation*}
 which is a contradiction.
 We conclude that $h(x,\alpha)=1$ on some interval in $(x_1,x_2)$.
 Since the function $e^t/t$ is strictly decreasing in $t\in(0,1)$, we have
 \begin{multline*}
  \int_{x_1}^{x_2} (u'(x)v(x)-u(x)v'(x))' dx \\
  = \lambda \int_{x_1}^{x_2} h(x,\alpha) 
     \left(\frac{e^{v(x)}}{v(x)} - \frac{e^{u(x)}}{u(x)} \right) u(x)v(x) dx
  > 0.
 \end{multline*}
 On the other hand, recalling $u(x_1)=v(x_1)\ge0$ and $u(x_2)=v(x_2)\ge0$, 
 we get
 \begin{align*}
  \int_{x_1}^{x_2} & (u'(x)v(x)-u(x)v'(x))' dx \\
    & = (u'(x_2)v(x_2) - u(x_2)v'(x_2)) - 
        (u'(x_1)v(x_1) - u(x_1)v'(x_1)) \\
    & = u(x_2)(u'(x_2) - v'(x_2)) - 
        u(x_1)(u'(x_1) - v'(x_1)) \le 0.
 \end{align*}
 This is a contradiction.
\end{proof}

\begin{lem}\label{nonexistence}
 Let $\lambda\in(0,\Lambda(1)]$.
 If $u$ is a positive non-even solution of \eqref{P}, then $||u||_\infty>1$.
\end{lem}

\begin{proof}
 Let $\beta\in(0,1]$ satisfy $\Lambda(\beta)=\lambda$.
 Since $U(x;\beta)$ is a positive even solution of 
 \eqref{P} and satisfies $||U(\,\cdot\,,\beta)||=\beta \le 1$.
 By Lemma \ref{||u||<1}, if \eqref{P} has a positive non-even solution $u$,
 then $u$ must satisfy $||u||_\infty>1$.
\end{proof}

We define $T(t,v)$ by
\[
 T(t,v)(x)
 = \int_{-1}^1 G(x,y) \Lambda(e^t) h(y,\alpha) e^{U(y;e^t)} (e^{v(y)}-1) dy,
\]
where $G(x,y)$ is a Green's function of
the operator $L[v]=-v''$ with $v(-1)=v(1)=0$, that is, 
\[
 G(x,y) =
  \left\{
   \begin{array}{ll}
    (1+x)(1-y)/2, & -1 \le x \le y \le 1, \\[1ex]
    (1-x)(1+y)/2, & -1 \le y \le x \le 1.
   \end{array}
  \right.
\]
We consider the equation
\begin{equation}
 v-T(t,v)=0.
  \label{v-T=0}
\end{equation}
Then  \eqref{v-T=0} has a trivial solution $v=0$ and
if $v$ is a solution of \eqref{v-T=0} at $t=\log\beta$ for some $\beta>0$, 
then $u(x):=U(x;\beta)+v(x)$ is a solution of \eqref{P} 
at $\lambda=\Lambda(\beta)$.

We will employ the following Rabinowitz's
global bifurcation theorem to $T(t,v)$.
See, for example, \cite[Theorem 2.5]{LS} and \cite[Theorem IV.12]{ST}. 

\begin{prop}\label{bifurcation}
 Let $E$ be a real Banach space and $T$: $\mathbb{R}\times E \to E$ compact 
 such that $T(t,0)=0$ for all $t \in \mathbb{R}$.
 Suppose that there exist constants $a$, $b \in \mathbb{R}$ with $a<b$ such
 that $(a,0)$ and $(b,0)$ are not bifurcation points for equation \eqref{v-T=0}.
 Furthermore, assume that
 \begin{equation*}
      \mbox{\rm deg}_{\rm LS}(I-T(a,\,\cdot\,),B_r(0),0)
  \ne \mbox{\rm deg}_{\rm LS}(I-T(b,\,\cdot\,),B_r(0),0),
 \end{equation*}
 where $B_r(0):=\{ v \in E : \| v \|_E < r \}$ is an isolating
 neighborhood of the trivial solution for both constants $a$ and $b$.
 Let
 \begin{equation*}
 {\mathcal S} = \overline{\{(t,v)
   : (t,v)\ \mbox{is\ a\ solution\ of\ \eqref{v-T=0}\ with} \ v \ne 0 \}}
   \cup ([a,b] \times \{0\})
 \end{equation*}
 and let ${\mathcal C}$ be the maximal connected subset of
 ${\mathcal S}$ containing $[a,b]\times \{0\}$.
 Then either
 \begin{enumerate}
  \item ${\mathcal C}$ is unbounded in $\mathbb{R} \times E$, or
  \item ${\mathcal C} \cap [(\mathbb{R}\setminus[a,b]) \times \{0\}] \ne \emptyset$.
 \end{enumerate}
\end{prop}

According to the condition of the change of the Leray-Schauder degree 
in Proposition \ref{bifurcation}, 
we find a bifurcation point in $[a,b]\times\{0\}$.

We will use Proposition \ref{bifurcation} with $E=C^1_0[-1,1]$.
(Hereafter, we regard $E$ as $C_0^1[-1,1]$.)
We easily see that $T$: $\mathbb{R}\times C_0^1[-1,1]\to C_0^1[-1,1]$ is 
a compact operator.
By the same arguments as in Kajikiya \cite[Section 5]{Kaj2023}, 
we obtain the following result.

\begin{lem}\label{deg}
 Let $t\in\mathbb{R}$. 
 Assume that $U(x;e^t)$ is nondegenerate.
 Then
 \begin{equation*}
  \mbox{\rm deg}_{\rm LS}(I-T(t,\,\cdot\,),B_r(0),0) = (-1)^{m(e^t)}
 \end{equation*}
 for each sufficiently small $r>0$.
\end{lem}

Lemmas \ref{mu1} and \ref{mu2} imply that, for each sufficiently small 
$\delta>0$, $U(x;\beta_2-\delta)$ and $U(x;\beta_2+\delta)$ are 
nondegenerate and 
\begin{equation*}
 m(\beta_2-\delta) = 1, \quad m(\beta_2+\delta)=2,
\end{equation*}
where $\beta_2$ is the constant determined in Lemma \ref{phi2}.
Now we fix such a $\delta>0$ and set $a=\log(\beta_2-\delta)$ and 
$b=\log(\beta_2+\delta)$.
Then $(a,0)$ and $(b,0)$ are not bifurcation points for equation \eqref{v-T=0}
and Lemma \ref{deg} means that
\begin{align*}
 & \mbox{\rm deg}_{\rm LS}(I-T(a,\,\cdot\,),B_r(0),0) = (-1)^1 = -1, \\ 
 & \mbox{\rm deg}_{\rm LS}(I-T(b,\,\cdot\,),B_r(0),0) = (-1)^2 = 1.
\end{align*}
From Proposition \ref{bifurcation}, it follows that (i) or (ii) of 
Proposition \ref{bifurcation} holds.

\begin{lem}\label{(beta1,0)}
 Let $\beta_2$ be as in Lemma \ref{phi2} and 
 let $t \in \mathbb{R}\setminus\{\log\beta_2\}$.
 Then $(t,0)$ is not a bifurcation point for equation \eqref{v-T=0}.
\end{lem}

\begin{proof}
 Lemmas \ref{mu1}, \ref{mu3} and \ref{mu2} imply that
 $(t,0)$ is not a bifurcation point for \eqref{v-T=0} 
 when $t\ne\log\beta_1$ or $t\ne\log\beta_2$.
 We will prove that $(\log\beta_1,0)$ is not a bifurcation point.

 Let $\phi_1$ be an eigenfunction of \eqref{Morse} with $\beta=\beta_1$ 
 corresponding to $\mu_1(\beta_1)$.
 We consider the map $F$: 
 $(0,\infty) \times (H_0^1(-1,1) \cap H^2(-1,1)) \to L^2(-1,1)$ defined by
 \begin{equation*}
   F(\lambda,u) = u'' + \lambda h(x,\alpha)e^u.
 \end{equation*}
 We check that $F$ at $(\Lambda(\beta_1),U(\,\cdot\,;\beta_1))$ satisfies 
 the assumptions of the Crandall-Rabinowitz bifurcation theorem 
 \cite[Theorem 3.2]{CR}.
 We observe that
 \begin{equation*}
  F_u(\Lambda(\beta_1),U(\,\cdot\,;\beta_1))w 
  = w'' + \Lambda(\beta_1) h(x,\alpha)e^{U(x;\beta_1)} w
 \end{equation*}
 for $w\in H_0^1(-1,1) \cap H^2(-1,1)$.
 Since $\mu_1(\beta_1)=0$, we conclude that the null-space 
 $N(F_u(\Lambda(\beta_1),U(\,\cdot\,;\beta_1)))$ can be expressed as 
 \begin{equation*}
  N(F_u(\Lambda(\beta_1),U(\,\cdot\,;\beta_1)))
  = \{ c\phi_1 : c \in \mathbb{R} \},
 \end{equation*}
 which means that $\dim N(F_u(\Lambda(\beta_1),U(\,\cdot\,;\beta_1)))=1$.

 Assume that $y \in R(F_u(\Lambda(\beta_1),U(\,\cdot\,;\beta_1)))$, that is, 
 there exists $w\in H_0^1(-1,1) \cap H^2(-1,1)$ such that
 \begin{equation*}
  w''(x) + \Lambda(\beta_1) h(x,\alpha)e^{U(x;\beta_1)} w(x) = y(x).
 \end{equation*}
 Then 
 \begin{equation*}
  R(F_u(\Lambda(\beta_1),U(\,\cdot\,;\beta_1))) 
   = \left\{ y \in L^2(-1,1) : 
       \int_{-1}^1 y(x) \phi_1(x) dx = 0 \right\}
 \end{equation*}
 and $\mbox{codim}\,R(F_u(\Lambda(\beta_1),U(\,\cdot\,;\beta_1)))=1$.
 We have
 \begin{equation*}
  F_\lambda(\Lambda(\beta_1),U(\,\cdot\,;\beta_1))
  = h(x,\alpha) e^{U(x;\beta_1)}.
 \end{equation*}
 Since $\phi_1$ has no zero in $(-1,1)$, we see that
 \begin{align*}
  \int_{-1}^1 h(x,\alpha) e^{U(x;\beta_1)} \phi_1(x) dx \ne 0,
 \end{align*}
 and hence
 \begin{equation*}
  F_\lambda(\Lambda(\beta_1),U(\,\cdot\,;\beta_1)) 
  \not\in R(F_u(\Lambda(\beta_1),U(\,\cdot\,;\beta_1))).
 \end{equation*}
 We note that $F(\Lambda(\beta_1),U(\,\cdot\,;\beta_1))=0$.
 From the Crandall-Rabinowitz bifurcation theorem 
 \cite[Theorem 3.2]{CR}, it follows that the solutions of 
 $F(\lambda,u)=0$ near $(\Lambda(\beta_1),U(\,\cdot\,;\beta_1))$ form a curve
 \begin{equation*}
  (\Lambda(\beta_1)+\tau(s),U(\,\cdot\,;\beta_1)+s \phi_1 +z(s)),
 \end{equation*}
 where $s\to(\tau(s),z(s))$ is a continuously differentiable function near 
 $s=0$ and $\tau(0)=\tau'(0)=z(0)=z'(0)=0$.
 On the other hand, $\{(\Lambda(\beta),U(\,\cdot\,;\beta)):\beta>0\}$ is 
 also a solution curve of $F(\lambda,u)=0$ through 
 $(\Lambda(\beta_1),U(\,\cdot\,;\beta_1))$.
 Hence, the solutions of \eqref{v-T=0} near $(\log\beta_1,0)$ form 
 $\{(t,0): |t-\log\beta_1|<\delta \}$ for some $\delta>0$.
 Consequently, $(\log\beta_1,0)$ is not a bifurcation point for \eqref{v-T=0}. 
\end{proof}

By Lemma \ref{(beta1,0)}, we conclude that (ii) of 
Proposition \ref{bifurcation} does not occur.
Therefore, there exists an unbounded connected set 
$\mathcal{C}'\subset \mathbb{R}\times (C^1_0[-1,1]\setminus \{0\})$ which 
emanates from $(\log\beta_2,0)$ such that each $(t,v)\in \mathcal{C}'$ is a 
solution of \eqref{v-T=0} with $v\ne0$.
We note that $(\log \beta_2,0)\not\in\mathcal{C}'$.

\begin{lem}\label{U+visPS}
 Let $(t,v) \in \mathcal{C}'$. 
 Then $U(x;e^t)+v(x)$ is a positive solution of \eqref{P} with 
 $\lambda=\Lambda(e^t)$.
\end{lem}

\begin{proof}
 We set $u(x)=U(x;e^t)+v(x)$.
 Then $u$ satisfies
 \begin{align*}
  u(x) -U(x;e^t) & = v(x) = T(t,v) = T(t,u-U(\,\cdot\,;e^t)) \\
   & = \int_{-1}^1 G(x,y) \Lambda(e^t) h(y,\alpha) e^{U(y;e^t)}
        (e^{u(y)-U(y;e^t)}-1) dy \\
   & = \int_{-1}^1 G(x,y) \Lambda(e^t) h(y,\alpha) e^{u(y)} dy \\ 
   & \quad \ -\int_{-1}^1 G(x,y) \Lambda(e^t) h(y,\alpha) e^{U(y;e^t)} dy \\
   & = \int_{-1}^1 G(x,y) \Lambda(e^t) h(y,\alpha) e^{u(y)} dy - U(x;e^t),   
 \end{align*}
 that is,
 \begin{equation*}
  u(x) = \int_{-1}^1 G(x,y) \Lambda(e^t) h(y,\alpha) e^{u(y)} dy.
 \end{equation*}
 This means that $u$ is a solution of \eqref{P} with $\lambda=\Lambda(e^t)$ 
 and $u(x)\ge 0$ for $x\in[-1,1]$.
 Since $u''(x)=0$ for $x \in (-\alpha,\alpha)$, $u''(x)<0$ for 
 $x \in (-1,-\alpha)\cup(\alpha,1)$ and $u(-1)=u(1)=0$, 
 we find that $u(x)>0$ for $x \in (-1,1)$.
\end{proof}

\begin{lem}\label{U+visPNES}
 Let $(t,v) \in \mathcal{C}'$. 
 Then $U(x;e^t)+v(x)$ is a positive non-even solution of \eqref{P} with 
 $\lambda=\Lambda(e^t)$.
\end{lem}

\begin{proof}
 Assume to the contrary that there exists $(t,v) \in \mathcal{C}'$
 such that $U(x;e^t)+v(x)$ is not a positive non-even solution of \eqref{P} 
 with $\lambda=\Lambda(e^t)$.
 Since $(\log \beta_2,0)\not\in\mathcal{C}'$, we have
 $(t,v)\ne(\log \beta_2,0)$.
 By Lemma \ref{U+visPS}, $U(x;e^t)+v(x)$ is a positive even solution of 
 \eqref{P} with $\lambda=\Lambda(e^t)$.
 Proposition \ref{evensolutions} implies that 
 $U(x;e^t)+v(x)=U(x;e^s)$, where $e^s=||U(\,\cdot\,;e^t)+v(\,\cdot\,)||_\infty$. 
 Since $v\ne0$, we have $t\ne s$.
 Therefore, $\mathcal{C}'$ has a connected component 
 $\mathcal{C}'' \subset \mathcal{C}'$ connecting $(\log \beta_2,0)$
 and some point $(s',0)\in\{(t,0): t\in\mathbb{R}\setminus\{\log \beta_2\}\}$, 
 that is, $(s',0)$ is a bifurcation point for equation \eqref{v-T=0}.
 This contradicts Lemma \ref{(beta1,0)}.
\end{proof}

Now we define the set $\mathcal{C}$ by
\begin{equation*}
 \mathcal{C}= \{ (\Lambda(e^t),U(\,\cdot\,;e^t)+v) : (t,v) \in \mathcal{C}' \}. 
\end{equation*}
Then $\mathcal{C}$ is a connected set which emanates from 
$(\Lambda(\beta_2),U(\,\cdot\,;\beta_2))$.
Lemmas \ref{lemlambda<} and \ref{U+visPNES} imply that 
if $(\lambda,u)\in \mathcal{C}$, then 
$u$ is a non-even positive solution of \eqref{P} and 
$\Lambda(||u||_\infty) \le \lambda < 4 \Lambda(||u||_\infty)$ holds.

Finally, we prove that there exist a sequence 
$\{(\lambda_n,u_n)\}\subset\mathcal{C}$ such that 
$\lambda_n\to0$ and $||u_n||_\infty\to\infty$.
Since $\mathcal{C}'$ is unbounded, there exists an unbounded sequence
$\{(t_n,v_n)\}_{n=1}^\infty \subset \mathcal{C}'$.
We set $u_n(x)=U(x;e^{t_n})+v_n(x)$ and $\lambda_n=\Lambda(e^{t_n})$.
Since $\Lambda(\beta)$ is a positive bounded function, 
by Lemmas \ref{lemlambda<} and \ref{U+visPNES}, 
there exists a constant $M>0$ such that $0<\lambda_n \le M$ for $n\ge 1$,
and hence $\{\lambda_n\}_{n=1}^\infty$ has a convergent subsequence.
We denote it by $\{\lambda_n\}_{n=1}^\infty$ again and its limit by 
$\lambda_0\in[0,M]$.
We claim that $\lambda_0=0$.
Assume $\lambda_0\in(0,M]$.
We have
\begin{equation*}
 \Lambda(||U(\,\cdot\,\,;e^{t_n})||_\infty) = \Lambda(e^{t_n}) = \lambda_n.
\end{equation*}
Lemma \ref{lemlambda<} yields
\begin{equation*}
 \Lambda(||u_n||_\infty) > \lambda_n/4.
\end{equation*}
Recalling \eqref{limlambda}, we conclude that
$\{t_n\}_{n=1}^\infty$,
$\{||U(\,\cdot\,;e^{t_n})||_\infty\}_{n=1}^\infty$ and
$\{||u_n||_\infty\}_{n=1}^\infty$ are bounded.
If $u$ is a positive solution of \eqref{P}, then
\begin{equation*}
 u'(x) = \int_m^x u''(y) dy 
         = - \lambda \int_m^x h(y,\alpha) e^{u(y)} dy, \quad x \in [-1,1],
\end{equation*}
where $m\in(-1,1)$ is a number with $u'(m)=0$,
and Lemma \ref{lemlambda<} shows
\begin{align*}
 |u'(x)| 
  & \le \lambda \int_{-1}^1 |h(y,\alpha)| e^{||u||_\infty} dy \\
  & < 8 \Lambda(||u||_\infty) e^{||u||_\infty}, \quad x \in [-1,1].
\end{align*}
Applying this inequality with $u(x)=U(x;e^{t_n})$ and $u(x)=u_n(x)$,
we conclude that $\{||U'(\,\cdot\,;e^{t_n})||_\infty\}_{n=1}^\infty$ and 
$\{||u_n'||\}_{n=1}^\infty$ are bounded, and hence
$\{||v_n'||_\infty\}_{n=1}^\infty$ is also bounded.
Consequently, $\{(t_n,v_n)\}_{n=1}^\infty \subset \mathcal{C}'$ is bounded.
This is a contradiction.
We have $\lambda_0=0$ as claimed.

By $\lambda_0=0$ and Lemma \ref{nonexistence}, we find that 
$||u_n||_\infty>1$ for all large $n$. 
Since $\Lambda(||u_n||_\infty)\le \lambda_n$ by Lemma \ref{lemlambda<},
from \eqref{limlambda}, it follows that $||u_n||_\infty\to\infty$ 
as $n\to\infty$.

The proof of Theorem \ref{main} is complete. 


\bigskip

{\bf Acknowledgments.} 
The authors would like to thank the reviewers for careful reading and 
helpful comments on important revisions that improve our manuscript.

\end{document}